\documentclass[11pt, twoside]{article}
\usepackage{mfpaperstuff}
\usepackage{needspace}
\usepackage{longtable}

\newgeometry{margin=0.421in}
\begin{document}

\newcommand{\thmcite}[2]{\textup{\textbf{\cite[#2]{#1}}}\ }

\newcommand\arrowover[1]{\scriptstyle\overrightarrow{\displaystyle #1}}
\newcommand\binoq[2]{\genfrac{[}{]}{0pt}{}{#1}{#2}}
\newcommand\mbinoq[2]{\Bigl[\medmath{\genfrac{}{}{-2pt}{}{#1}{#2}}\Bigr]}
\newcommand{\ol}{\overline}
\newcommand{\ul}{\tilde}
\newcommand{\sss}{W_}
\newcommand{\nchar}{\operatorname{char}}
\newcommand{\hhh}{\scrh_}
\renewcommand\hom{\operatorname{Hom}}
\newcommand\spe{c_\xi}
\newlength\adju
\newlength\bdju
\setlength\bdju{5pt}
\newcommand\tsum[1]{\settowidth{\adju}{\hbox{\scriptsize $#1$}}\addtolength\adju{-\bdju}\setlength\adju{0.45\adju}\hspace{-0.8\adju}\sum_{w\in#1}\hspace{-\adju}T_w}
\newcommand\ttsum[1]{\sum_{w\in#1}\!\!T_w}
\newcommand\lsum[1]{\settowidth{\adju}{\hbox{\scriptsize $#1$}}\addtolength\adju{-\bdju}\setlength\adju{0.45\adju}\hspace{-0.8\adju}\sum_{w\in#1}\hspace{-\adju}q^{l(w)}}
\newcommand\sh[1]{\Gamma_{#1}}
\newcommand\shif[2]{\sh{#1}(#2)}
\newcommand\shiflr[2]{\sh{#1}\left(#2\right)}
\newcommand\rf[1]{|#1\rangle}
\newcommand\lf[1]{\langle#1|}
\newcommand\ptn\vdash
\newcommand\compn\vDash
\newcommand\caltr[2]{\calt_{\hspace{-2pt}\operatorname{r}}(#1,#2)}
\newcommand\calto[2]{\calt_{\hspace{-2pt}0}(#1,#2)}
\newcommand\hone{(H1)}
\newcommand\htwo{(H2)}
\newcommand\hthree{(H3)}
\newcommand\hfour{(H4)}
\newcommand\rowq{^\diamond}

\pdfpagewidth=12.5cm
\pdfpageheight=8.839cm
\setcounter{page}0
{\footnotesize This is the author's version of a work that was accepted for publication in the Journal of Algebra. Changes resulting from the publishing process, such as peer review, editing, corrections, structural formatting, and other quality control mechanisms may not be reflected in this document. Changes may have been made to this work since it was submitted for publication. A definitive version was subsequently published in\\
\textit{J.\ Algebra} \textup{\textbf{364} (2012) 38--516.\\
http://dx.doi.org/10.1016/j.jalgebra.2012.04.020}\normalsize}
\newgeometry{margin=1in,includehead,includefoot}

\title{An algorithm for semistandardising homomorphisms}
\msc{20C30, 20C08, 05E10}
\toptitle

\begin{abstract}
Suppose $\mu$ is a partition of $n$ and $\la$ a composition of $n$, and let $S^\mu$, $M^\la$ denote the Specht module and permutation module defined by Dipper and James for the Iwahori--Hecke algebra $\hhh n$ of the symmetric group $\sss n$.  We give an explicit fast algorithm for expressing a tableau homomorphism $\hat\phi_A:S^\mu\to M^\la$ as a linear combination of semistandard homomorphisms.  Along the way we provide a utility result related to removing rows from tableaux.
\end{abstract}

\pdfpagewidth=8.27in
\pdfpageheight=11.69in

\renewcommand\baselinestretch{1.172}

\Yboxdim{12pt}
\section{Introduction}

Throughout this paper, $\bbf$ is an arbitrary field and $q$ is a fixed non-zero element of $\bbf$.  Given a positive integer $n$, $\sss n$ denotes the symmetric group of degree $n$ and $\hhh n$ denotes the Iwahori--Hecke algebra of $\sss n$: this is the unital associative $\bbf$-algebra with generators $T_1,\dots,T_{n-1}$ and relations
\begin{alignat*}2
(T_i-q)(T_i+1)&=0&&\text{for }1\ls i\ls n{-}1,\\
T_iT_j&=T_jT_i&&\text{for }1\ls i\ls j{-}2\ls n{-}3,\\
T_iT_{i+1}T_i&=T_{i+1}T_iT_{i+1}&\qquad&\text{for }1\ls i\ls n{-}2.
\end{alignat*}
The representation theory of $\hhh n$ closely resembles (and informs) the modular representation theory of $\sss n$.  This theory was first explored in detail by Dipper and James, whose account \cite{dj} we closely follow.  In particular, they define a `permutation module' $M^\la$ for each composition $\la$ of $n$, and a Specht module $S^\mu$ for each partition $\mu$ of $n$.  Understanding the structure of Specht modules is an important goal in understanding the representation theory of $\hhh n$, and one well-studied aspect of this is the examination of $\hhh n$-homomorphisms between Specht modules.

The aim of this paper is to provide a `straightening' result which yields (as long as $q\neq-1$) a fast algorithm to determine $\hom_{\hhh n}(S^\mu,S^\la)$.  This result actually concerns the space of homomorphisms from $S^\mu$ to $M^\la$; an important family of such homomorphisms is labelled by $\mu$-tableaux of type $\la$; when $q\neq-1$, this family spans the space $\hom_{\hhh n}(S^\mu,M^\la)$ \cite[Corollary 8.2]{dj2}. Regardless of $q$, the space spanned by the tableau homomorphisms has a basis consisting of the homomorphisms labelled by \emph{semistandard} tableaux. However, existing proofs of this result involve a loss of information which means that it is not easy to express a given tableau homomorphism explicitly in terms of semistandard homomorphisms.  In this paper we prove a relation between tableau homomorphisms which allows us to do exactly this; this closely resembles the \emph{Garnir relations}  which are typically used to straighten elements of the Specht module.

In the next section, we very briefly recall the background information that we need, referring to \cite{dj} for most of it.  In Section \ref{s:twopart} we prove our Garnir-like relation for tableaux with only two rows.  In Section \ref{s:rowrem} we prove the necessary results concerning row removal to show that our Garnir relation holds generally.  Finally, in Section \ref{s:algo} we give an algorithm for expressing a tableau homomorphism as a linear combination of semistandard homomorphisms.

\section{Background}

\subs{Basic definitions}

We take almost all of our notation from \cite{dj}, with only very minor modifications.  We provide a brief index here of most of the notation, referring the reader to \cite{dj} for definitions.  Homomorphisms are considered in more detail below.  Note that we follow \cite{dj} by considering right modules and having all functions act on the right.

\begin{longtable}{ll}
$\sss n$&the symmetric group on $\{1,\dots,n\}$\\
$l$&the (Coxeter) length function on $\sss n$\\
$\la\compn n$&$\la$ is a composition of $n$\\
$\la\ptn n$&$\la$ is a partition of $n$\\
$\la'$&the conjugate partition to $\la$\\
$W_\la$&the standard Young subgroup of $\sss n$ defined by $\la\compn n$\\
$\scrd_\la$&the set of minimal-length right coset representatives for $W_\la$ in $\sss n$\\
$\sim$&the row-equivalence relation on tableaux\\
$\fkt^\la$&the $\la$-tableau which has $1,\dots,n$ in order along successive rows, for $\la\compn n$\\
$\fkt_\la$&the $\la$-tableau which has $1,\dots,n$ in order down successive columns, for $\la\ptn n$\footnote{n.b.\ this notation is not defined in \cite{dj}, although this tableau is used there.  This notation is standard elsewhere.}\\
$w_\la$&the permutation such that $\fkt^\la w_\la=\fkt_\la$\\
$\hhh n$&the Iwahori--Hecke algebra of $\sss n$ over $\bbf$ with parameter $q$\\
$T_1,\dots,T_{n-1}$&the standard generators of $\hhh n$\\
$T_w$&standard basis element of $\hhh n$, for $w\in\sss n$\\
$x_\la$&$\ttsum{W_\la}$\\
$y_\la$&$\sum_{w\in W_\la}(-q)^{-l(w)}T_w$\\
$M^\la$&the `permutation module' $x_\la\hhh n$\\
$S^\la$& the Specht module $x_\la T_{w_\la}y_{\la'}\hhh n$
\end{longtable}

\subs{Multisets}

In this paper we frequently employ multisets of positive integers.  If $C$ is a multiset of positive integers, then we write $C_i$ for the number of $i$s in $C$.  Given multisets $C,D$, we write $C\sqcup D$ for the multiset with $(C\sqcup D)_i=C_i+D_i$ for each $i$.  If $A$ is a tableau, we write $A^i$ for the multiset of entries in row $i$ of $A$.

\subs{Quantum binomial coefficients}

Given our fixed $q$, we write $[n]=1+q+\dots+q^{n-1}$ for any non-negative integer $n$, and $[n]!=\prod_{i=1}^n[i]$.  These expressions are useful in the study of $\hhh n$, since we have $[n]!=\sum_{w\in\sss n}q^{l(w)}$.  Furthermore, if $\la\compn n$, then any $w\in\sss n$ can be written uniquely as $vd$ with $v\in W_\la$, $d\in\scrd_\la$ and $l(w)=l(v)+l(d)$, so we have
\[
\lsum{\scrd_\la}=\frac{[n]!}{\prod_i[\la_i]!}.
\]
Note that the right-hand side of this expression is a polynomial in $q$, so makes sense even when $\prod_i[\la_i]!$ is zero.

If $q$ is an indeterminate and $n\gs r\gs 0$, then we can define the \emph{quantum binomial coefficient}
\[
\binoq nr=\frac{[n]!}{[r]![n-r]!}.
\]
In fact, this is a polynomial in $q$, so we can extend the definition to arbitrary $q$ by defining $\binoq nr$ to be the specialisation of this polynomial.

\subs{Homomorphisms between permutation modules}

We now recall the definitions from \cite[\S3]{dj} on homomorphisms between permutation modules. Throughout this subsection we fix $\la,\mu\compn n$.

Recall that a $\mu$-tableau $A$ of type $\la$ is \emph{row-standard} if its entries increase weakly along the rows. We write $\caltr\mu\la$ for the set of row-standard $\mu$-tableaux of type $\la$.  If in addition $\mu\ptn n$, then we say that a $\mu$-tableau of type $\la$ is \emph{semistandard} if its entries are weakly increasing along the rows and strictly increasing down the columns, and we write $\calto\mu\la$ for the set of semistandard $\mu$-tableaux of type $\la$.

Given a $\mu$-tableau $A$ of type $\la$, we define the permutation $1_A$ by specifying that $\fkt^\la1_A$ is the row-standard $\la$-tableau (of type $1^n$) in which $i$ belongs to row $r$ if the position occupied by $i$ in $\fkt^\mu$ is occupied by $r$ in $A$.  For example, if
\Yvcentermath1
\begin{align*}
A&=\young(123,112)\\
\intertext{then}
\fkt^\la1_A&=\young(145,26,3)
\end{align*}
so that $1_A=(2\ 4)(3\ 5\ 6)$.

The map $A\mapsto1_A$ defines a bijection from $\caltr\mu\la$ to the set $\scrd_\la\cap\scrd^{-1}_\mu$ of minimal-length ($W_\la$,$W_\mu$) double coset representatives in $\sss n$. For later use we record the simple fact that if $A\in\caltr\mu\la$, then 
\[
l(1_A)=\sum_{g<h}\sum_{i<j}A^g_jA^h_i.
\]

Now we consider homomorphisms. If $A\in\caltr\mu\la$, then we define an $\hhh n$-homomorphism $\phi_A:M^\mu\to M^\la$ by
\[
x_\mu\phi_A=\tsum{W_\la1_AW_\mu}.\tag*{\hone}\\
\]
We call this a \emph{tableau homomorphism}.  In the case where $\mu$ is a partition, we write $\hat\phi_A$ for the restriction of $\phi_A$ to the Specht module $S^\mu$. We remark that in more recent literature $\phi_A$, $\hat\phi_A$ are written as $\Theta_A$, $\hat\Theta_A$.

Now we have the following theorem.

\needspace{0pt}\begin{thm}\textup{\textbf{\cite[Theorem 3.4]{dj}, \cite[Corollary 8.7]{dj2}}}\ \label{semi}\indent
\begin{enumerate}
\vspace{-\topsep}
\nopagebreak[4]
\item
If $\mu,\la\compn n$, then the set $\lset{\phi_A}{A\in\caltr\mu\la}$ is a basis for $\hom_{\hhh n}(M^\mu,M^\la)$.
\item
If $\mu\ptn n$ and $\la\compn n$, then the set $\lset{\hat\phi_A}{A\in\calto\mu\la}$ is a basis for the subspace of $\hom_{\hhh n}(S^\mu,M^\la)$ spanned by $\lset{\hat\phi_A}{A\in\caltr\mu\la}$.
\item
If $\mu\ptn n$, $\la\compn n$ and $q\neq-1$, then every $\hhh n$-homomorphism $S^\mu\to M^\la$ can be extended to a homomorphism $M^\mu\to M^\la$.  Hence $\lset{\hat\phi_A}{A\in\calto\mu\la}$ is a basis for $\hom_{\hhh n}(S^\mu,M^\la)$.
\end{enumerate}
\end{thm}

The focus of this paper is part (2) of the above theorem; we provide an explicit algorithm for writing a given homomorphism $\phi_A$ as a linear combination of semistandard homomorphisms. Of course, the proof of Theorem \ref{semi}(2) does yield an algorithm for `semistandardising' a tableau homomorphism, but it is rather slow, even in the symmetric group case $q=1$; one must consider the image of a polytabloid under the different homomorphisms $\phi_A$, and keep track of the coefficients of various tabloids in these images.  Our algorithm involves just manipulation of row-standard tableaux, forgetting the underlying homomorphisms. In simple computer experiments, our algorithm seems considerably quicker.

\smallskip

There are several other useful ways to write the image $x_\mu\phi_A$.  From \cite[p.~30]{dj}, we have
\[
x_\mu\phi_A=\sum_{A\rowq\sim A}x_\la T_{1_{A\rowq}}.\tag*{\htwo}
\]
By considering the row permutations sending $A$ to the various $A\rowq$, we can also write this (following \cite[Theorem 3.4]{dj}) as
\[
x_\mu\phi_A=x_\la T_{1_A}\tsum{\scrd_\nu\cap W_\mu},\tag*{\hthree}
\]
where $\nu$ is the `composition defined by reading $A$ along rows', i.e.\ the composition
\[
(A^1_1,A^1_2,A^1_3,\dots,A^2_1,A^2_2,A^2_3,\dots,A^3_1,A^3_2,\dots).
\]
We will often use the expression \hthree, and write the term $\ttsum{\scrd_\nu\cap W_\mu}$ as $\rf A$.

Symmetrically, we have
\[
x_\mu\phi_A=\tsum{\scrd^{-1}_\pi\cap W_\la}T_{1_A}x_\mu\tag*{\hfour},
\]
where $\pi$ is the composition
\[
(A^1_1,A^2_1,A^3_1,\dots,A^1_2,A^2_2,A^3_2,\dots,A^1_3,A^2_3,\dots).
\]
We shall write the factor $\ttsum{\scrd^{-1}_\pi\cap W_\la}$ as $\lf A$.

Our main interest is in understanding linear combinations $\sum_{i=1}^rc_i\phi_{A_i}$ (where $A_1,\dots,A_r\in\caltr\mu\la$ and $c_1,\dots,c_r\in\bbf$) for which $\sum_{i=1}^rc_i\hat\phi_{A_i}=0$.  Since $S^\mu$ is generated by $x_\mu T_{w_\mu}y_{\mu'}$, this condition amounts to saying that $\sum_{i=1}^r c_ix_\mu\phi_{A_i}$ is annihilated by $T_{w_\mu}y_{\mu'}$.

\section{A Garnir-type relation for tableaux with two rows}\label{s:twopart}

Now we come to our main result.  We begin with tableaux having only two rows.

\begin{thm}\label{main}
Suppose $\mu=(m,n-m)$ is a partition, and $R,S,T$ are multisets of positive integers with $|R|+|S|+|T|=n$ and $|S|>m$.  Let $\cals$ be the set of all pairs $(U,V)$ such that $S=U\sqcup V$ and $|U|=m-|R|$.  For each $(U,V)\in\cals$, let $A[U,V]$ denote the row-standard $\mu$-tableau with
\[
A[U,V]^1=R\sqcup U,\qquad A[U,V]^2=T\sqcup V.
\]
Then
\[
\sum_{(U,V)\in\cals}\prod_i\mbinoq{R_i+U_i}{R_i}\mbinoq{T_i+V_i}{T_i}\prod_{i<j}q^{R_jU_i+T_iV_j}\hat\phi_{A[U,V]}=0.
\]
\end{thm}

An example of this theorem is given in Section \ref{s:algo}.  The way we prove the theorem will be to express the given homomorphism (extended to $M^\mu$) as a composition $\psi\chi$ such that $\psi|_{S^\mu}=0$.  To do this, we need a few simple results concerning composition of tableau homomorphisms.  Some of these are probably known, but we prove everything in order to keep our account self-contained.

Our first result shows how to use a tableau of the form $\gyoung(;1;1_2\cdots;1;2_2\cdots;2)$ to merge two rows of a tableau.

\begin{propn}\label{bigpart}
Suppose $\xi=(r,m-r)$ is a two-part composition of $m$, and $C$ is a row-standard $\xi$-tableau of type $\alpha\compn m$.  Let $B$ be the unique row-standard $(m)$-tableau of type $\xi$, and $A$ the unique row-standard $(m)$-tableau of type $\alpha$.  Then
\[
\phi_B\phi_C = \prod_i\mbinoq{A^1_i}{C^1_i}\prod_{i<j}q^{C^1_jC^2_i}\phi_A.
\]
\end{propn}

\begin{pf}
Equation \hone{} shows that for a row-standard $(m)$-tableau $D$ of any type, we have
\[
x_{(m)}\phi_D=x_{(m)}.
\]
So we just need to show that $x_{(m)}\phi_B\phi_C$ is the appropriate scalar multiple of $x_{(m)}$. We have
\[
x_{(m)}\phi_B=x_\xi\tsum{\mathscr{D}_\xi}
\]
by \hthree, while \hfour{} gives
\[
x_\xi\phi_C=\tsum{\scrd_\pi^{-1}\cap W_\alpha}T_{1_C}x_\xi,
\]
where $\pi$ is the composition $(C^1_1,C^2_1,C^1_2,C^2_2,C^1_3,C^2_3,\dots)$. So
\[
x_{(m)}\phi_B\phi_C=\tsum{\scrd_\pi^{-1}\cap W_\alpha}T_{1_C}x_\xi\tsum{\mathscr{D}_\xi}=\tsum{\scrd_\pi^{-1}\cap W_\alpha}T_{1_C}x_{(m)}.
\]
To evaluate this expression, we use the fact that $T_wx_{(m)}=q^{l(w)}x_{(m)}$ for $w\in\sss n$.  The length of $1_C$ is $\sum_{i<j}C^1_jC^2_i$, and since $W_\pi\ls W_\alpha$, we have
\[
\lsum{\scrd_\pi^{-1}\cap W_\alpha}=\frac{\sum_{s\in W_\alpha}q^{l(s)}}{\sum_{s\in W_\pi}q^{l(s)}}=\frac{\prod_i[\alpha_i]!}{\prod_i[\pi_i]!}=\prod_i\mbinoq{A^1_i}{C^1_i}.
\]
The result follows.
\end{pf}

The next result is a more complicated version of the result above.

\begin{propn}\label{moresoph}
Suppose $(r,u,v,t)$ is a composition of $n$.  Let $C$ be a row-standard $(r,u,v,t)$-tableau of type $\la\compn n$, and let $B$ be the row-standard $(r+u,v+t)$-tableau in which there are $r$ $1$s and $u$ $2$s in the first row, and $v$ $3$s and $t$ $4$s in the second.  Let $A$ be the row-standard $(r+u,v+t)$-tableau with $A^1=C^1\sqcup C^2$ and $A^2=C^3\sqcup C^4$.  Then
\[
\phi_B\phi_C = \prod_i\mbinoq{A^1_i}{C^1_i}\mbinoq{A^2_i}{C^3_i}\prod_{i<j}q^{C^1_jC^2_i+C^3_jC^4_i}\phi_A.
\]
\end{propn}

\begin{pf}
Since $1_B=1$, \hthree{} gives
\[
x_{(r+u,v+t)}\phi_B=x_{(r,u,v,t)}\tsum{\scrd_{(r,u,v,t)}\cap W_{(r+u,v+t)}}
\]
and \hfour{} gives
\[
x_{(r,u,v,t)}\phi_C=\tsum{\scrd^{-1}_\pi\cap W_\la}T_{1_C}x_{(r,u,v,t)},
\]
where $\pi$ is the composition $(C^1_1,C^2_1,C^3_1,C^4_1,C^1_2,C^2_2,C^3_2,C^4_2,C^1_3,C^2_3,C^3_3,C^4_3,\dots)$. So
\begin{align*}
x_{(r,u,v,t)}\phi_B\phi_C&=\tsum{\scrd^{-1}_\pi\cap W_\la}T_{1_C}x_{(r,u,v,t)}\tsum{\scrd_{(r,u,v,t)}\cap W_{(r+u,v+t)}}.
\\
&=\tsum{\scrd^{-1}_\pi\cap W_\la}T_{1_C}x_{(r+u,v+t)}.
\end{align*}
$T_{1_C}$ factorises as $T_{1_A}T_d$, where $d\in W_{(r+u,v+t)}$, so the above expression becomes
\[
q^{l(d)}\tsum{\scrd^{-1}_\pi\cap W_\la}T_{1_A}x_{(r+u,v+t)}=q^{l(d)}\tsum{\scrd^{-1}_\rho\cap W_\la}\ \ \ \ \tsum{\scrd^{-1}_\pi\cap W_\rho}T_{1_A}x_{(r+u,v+t)}
\]
where $\rho$ is the composition $(A^1_1,A^2_1,A^1_2,A^2_2,\dots)$, and $l(d)=\sum_{i<j}(C^1_jC^2_i+C^3_jC^4_i)$.  Now $1_A$ conjugates the pair $(W_\pi,W_\rho)$ to the pair $(W_\sigma,W_\tau)$, where
\begin{align*}
\sigma&=(C^1_1,C^2_1,C^1_2,C^2_2,C^1_3,C^2_3,\dots,C^3_1,C^4_1,C^3_2,C^4_2,C^3_3,C^4_3,\dots),\\
\tau&=(A^1_1,A^1_2,A^1_3,\dots,A^2_1,A^2_2,A^2_3,\dots),
\end{align*}
and this gives
\[
\tsum{\scrd^{-1}_\pi\cap W_\rho}T_{1_A}=T_{1_A}\tsum{\scrd^{-1}_\sigma\cap W_\tau}.
\]
$W_\tau$ is contained in $W_{(r+u,v+t)}$, so
\[
\tsum{\scrd^{-1}_\sigma\cap W_\tau}x_{(r+u,v+t)}=\lsum{\scrd^{-1}_\sigma\cap W_\tau}x_{(r+u,v+t)}=\prod_i\mbinoq{A^1_i}{C^1_i}\mbinoq{A^2_i}{C^3_i}x_{(r+u,v+t)}.
\]
So we have
\begin{align*}
x_{(r,u,v,t)}\phi_B\phi_C&=\prod_i\mbinoq{A^1_i}{C^1_i}\mbinoq{A^2_i}{C^3_i}\prod_{i<j}q^{C^1_jC^2_i+C^3_jC^4_i}\tsum{\scrd^{-1}_\rho\cap W_\la}T_{1_A}x_{(r+u,v+t)}\\
&=\prod_i\mbinoq{A^1_i}{C^1_i}\mbinoq{A^2_i}{C^3_i}\prod_{i<j}q^{C^1_jC^2_i+C^3_jC^4_i}x_{(r,u,v,t)}\phi_A.\qhere{\sum_{i<j}}
\end{align*}
\end{pf}

The next result effectively enables us to split a row of a tableau in two, using a tableau of the form
\[
\gyoung(;1;1_3\cdots;1;1,;2;2_4\cdots;2;2,;2;2_3\cdots;2;2,;3;3_2\cdots;3;3).
\]

\begin{lemma}\label{easy}
Suppose $(r,u,v,t)$ is a composition of $n$, and let $D$ be the unique $(r,u,v,t)$-tableau in which the entries in row $1$ are equal to $1$, the entries in rows $2$ and $3$ are all equal to $2$, and the entries in row $4$ are all equal to $3$. Let $E$ be a row-standard $(r,u+v,t)$-tableau of any type, and let $\calc$ be the set of all row-standard $(r,u,v,t)$-tableaux $C$ such that
\[
C^1=E^1,\qquad C^2\sqcup C^3=E^2,\qquad C^4=E^3.
\]
Then
\[
\phi_D\phi_E=\sum_{C\in\calc}\phi_C.
\]
\end{lemma}

\begin{pf}
We have $x_{(r,u,v,t)}\phi_D=x_{(r,u+v,t)}$, so by \htwo
\[
x_{(r,u,v,t)}\phi_D\phi_E=x_\la\sum_{E\rowq\sim E}T_{1_{E\rowq}},
\]
where $\la$ is the type of $E$. But it is easy to see that
\[
\lset{1_{E\rowq}}{E\rowq\sim E} = \bigsqcup_{C\in\calc}\lset{1_{C\rowq}}{C\rowq\sim C}
\]
so
\[
x_{(r,u,v,t)}\phi_D\phi_E=x_\la\sum_{C\in\calc}\sum_{C\rowq\sim C}T_{1_{C\rowq}}=\sum_{C\in\calc}x_\la\sum_{C\rowq\sim C}T_{1_{C\rowq}}=\sum_{C\in\calc}x_\la\phi_C.\qhere{\sum_{C\in\calc}}
\]
\end{pf}

The preceding results allow us to express the homomorphism of Theorem \ref{main} as a composition of homomorphisms indexed by row-standard tableaux.

\begin{propn}\label{compo}
Suppose $R,T,S,\cals$ are as in Theorem \ref{main}, and let $r=|R|$, $s=|S|$, $t=|T|$, $u=m-r$, $v=n-m-t$. Let $E$ be the row-standard $(r,s,t)$-tableau which has $(E^1,E^2,E^3)=(R,S,T)$, let $B$ be the $(r+u,v+t)$-tableau from Proposition \ref{moresoph}, and let $D$ be the $(r,u,v,t)$-tableau from Lemma \ref{easy}.  Then
\[
\phi_B\phi_D\phi_E=\sum_{(U,V)\in\cals}\prod_i\mbinoq{R_i+U_i}{R_i}\mbinoq{T_i+V_i}{T_i}\prod_{i<j}q^{R_jU_i+T_iV_j}\phi_{A[U,V]}.
\]
\end{propn}

\begin{pf}
From Lemma \ref{easy}, we have
\[
\phi_D\phi_E=\sum_{C\in\calc}\phi_C;
\]
each tableau $C\in\calc$ satisfies $(C^1,C^2,C^3,C^4)=(R,U,V,T)$ for some pair $(U,V)\in\cals$.  When we compose $\phi_B$ with $\phi_C$ using Proposition \ref{moresoph}, the tableau $A$ we obtain is precisely $A[U,V]$, and the coefficient of $\phi_A$ is the coefficient of $\phi_{A[U,V]}$ in the proposition.
\end{pf}

\begin{pf}[Proof of Theorem \ref{main}]
Using Proposition \ref{compo}, we must show that the composition $\phi_B\phi_D\phi_E$ kills $S^\mu$.  In fact, we will show that $\phi_B\phi_D$ kills $S^\mu$.  The codomain of $\phi_B\phi_D$ is the permutation module $M^{(r,s,t)}$, while the Specht module $S^\mu$ is generated by $x_\mu T_{w_\mu}y_{\mu'}$, so it suffices to show that $M^{(r,s,t)}y_{\mu'}=0$.  But this is straightforward: $M^{(r,s,t)}$ is spanned by elements $x_{(r,s,t)}T_d$ which correspond to $(r,s,t)$-tabloids (this correspondence is introduced on \cite[p.\ 31]{dj}).  If $u$ is an $(r,s,t)$-tabloid, then because $s>m$ there must be some pair $(i,i{+}1)$ which lie in the same column of $\fkt_\mu$ and both lie in the second row of $u$.  $y_{\mu'}$ may be factorised as $(1-q^{-1}T_i)h$, and the fact that $i,i{+}1$ both lie in the same row of $u$ means that $uT_i=qu$.  Hence $uy_{\mu'}=u(1-q^{-1}T_i)h=0$.
\end{pf}

\begin{rmk}
The argument used here to show that $M^{(r,s,t)}y_{\mu'}=0$ is a special case of \cite[Lemma 4.1]{dj} which shows that $M^\nu y_{\mu'}=0$ whenever $\nu'\ndom\mu'$.
\end{rmk}

\section{Row removal}\label{s:rowrem}

In this section, we prove a result which allows us to apply our Garnir-type relation to tableaux with more than two rows.  Using the results of this section, we will be able to prove a $q$-analogue of \cite[Lemma 4]{fm}, which essentially says that given a linear relation between homomorphisms $\hat\phi_A$, we can add a fixed combination of rows to each of the tableaux and preserve the relation.

We remark that row-removal for tableau homomorphisms has been considered at length, by the author and Lyle \cite{fl} for the symmetric group, and by Lyle and Mathas \cite{LM} for the Hecke algebra.  However, the result we prove here seems not to have appeared before.

We need to introduce some more notation.  Given $m\ls n$ and $0\ls r\ls n-m$, we define the homomorphism $\sh r:\hhh m\to \hhh n$ by $T_i\mapsto T_{i+r}$.  So $\sh0$ is the usual embedding of $\hhh m$ in $\hhh n$.  Clearly we have $\sh r\sh s=\sh{r+s}$ for any $r,s$.  We will also occasionally use $\sh r$ to denote the corresponding map $\sss m\to\sss n$.

\subs{Adding a row at the top}

The main result we want to prove here is the following.

\begin{propn}\label{topmain}
Suppose $\mu\ptn n$ and $\la\compn n$, and $A_1,\dots,A_r$ are row-standard $\mu$-tableaux of type $\la$ which are identical in row $1$.  Write $\tilde\mu=(\mu_2,\mu_3,\dots)$, and let $\tilde A_1,\dots,\tilde A_r$ be the $\tilde\mu$-tableaux obtained by deleting the first row from each of $A_1,\dots,A_r$.  If $c_1,\dots,c_r\in\bbf$ are such that $\sum_{i=1}^rc_i\hat\phi_{\tilde A_i}=0$, then $\sum_{i=1}^rc_i\hat\phi_{A_i}=0$.
\end{propn}

In order to prove Proposition \ref{topmain}, we make two comparisons:
\begin{itemize}
\item
we compare $x_\mu\phi_A$ with $x_{\tilde\mu}\phi_{\tilde A}$, for a row-standard $\mu$-tableau $A$;
\item
we compare the term $T_{w_\mu}y_{\mu'}$ which defines the Specht module with $T_{w_{\tilde\mu}}y_{\tilde\mu'}$.
\end{itemize}

\begin{propn}\label{rowtophom}
Suppose $A$ is a row-standard $\mu$-tableau of type $\la$.  Write $\tilde\mu=(\mu_2,\mu_3,\dots)$, and let $\tilde A$ be the $\tilde\mu$-tableau obtained by deleting the first row of $A$.  Let $\tilde\la$ be the type of $\tilde A$.  Then we have
\[
x_\mu\phi_A=h\shif{\mu_1}{x_{\tilde\mu}\phi_{\tilde A}}
\]
for some $h\in\hhh n$ depending only on $\mu,\la,\tilde\la$.
\end{propn}

The way we prove Proposition \ref{rowtophom} is to express $\phi_A$ as a composition of tableau homomorphisms.  Let $B$ be the $\mu$-tableau which is obtained from $A$ by replacing every entry in the first row with a $1$ and adding $1$ to every entry below the first row, and let $C$ be the tableau whose first row is the same as the first row of $A$, and whose $i$th row consists of $\tilde\la_{i-1}$ $(i{-}1)$s, for $i\gs2$.  Note that both the type of $B$ and the shape of $C$ equal the composition $\kappa=(\mu_1,\tilde\la_1,\tilde\la_2,\dots)$.

\Yvcentermath1
\begin{eg}
Suppose $\mu=(4,4,2,1)$, $\la=(3,4,2,1,1)$ and
\[
A=\young(1233,1245,22,1).
\]
Then $\tilde\la=(2,3,0,1,1)$, and
\[
B=\young(1111,2356,33,2),\qquad C=\gyoung(;1;2;3;3,;1;1,;2;2;2,:-,;4,;5).
\]

\end{eg}

\begin{lemma}\label{comptop}
With the notation above, we have $\phi_B\phi_C=\phi_A$.
\end{lemma}

\begin{pf}
Using \hthree, we get
\[
x_\mu\phi_B\phi_C=x_\la T_{1_C}\rf CT_{1_B}\rf B.
\]
Since the entries in the $i$th row of $C$ are constant for $i\gs2$, we have $\rf C\in\hhh{\mu_1}$.  On the other hand, since the entries in the first row of $B$ all equal $1$, we have
\[
T_{1_B}=\shif{\mu_1}{T_{1_{\tilde A}}}\in\shif{\mu_1}{\hhh{n-\mu_1}}.
\]
So $\rf C$ and $T_{1_B}$ commute.  It is not too hard to see that $T_{1_C}T_{1_B}=T_{1_A}$, and it is immediate from the definition of $\rf A$ that $\rf C\rf B=\rf A$, so
\[
x_\mu\phi_B\phi_C=x_\la T_{1_A}\rf A=x_\mu\phi_A.\qhere{x_\mu}
\]
\end{pf}

\begin{pf}[Proof of Proposition \ref{rowtophom}]
We can write $x_\kappa\phi_C$ as $\lf CT_{1_C}x_\kappa$, so
\[
x_\mu\phi_A=x_\mu\phi_B\phi_C=\lf CT_{1_C}x_\kappa T_{1_B}\rf B.
\]
We have $x_\kappa=x_{(\mu_1)}\shif{\mu_1}{x_{\tilde\la}}$; as noted in the proof of Lemma \ref{comptop}, $T_{1_B}$ is just $\shif{\mu_1}{T_{1_{\tilde A}}}$, and it is immediate that $\rf B=\shif{\mu_1}{\rf{\tilde A}}$.  Putting this together, we get
\[
x_\mu\phi_A=\lf CT_{1_C}x_{(\mu_1)}\shiflr{\mu_1}{x_{\tilde\la}T_{1_{\tilde A}}\rf{\tilde A}};
\]
the tableau $C$ depends only on $\mu,\la,\tilde\la$, so the result follows.
\end{pf}

The next thing we need in order to prove Proposition \ref{topmain} is the following.

\begin{propn}\label{tytop}
Suppose $\mu$ is a partition of $n$, and write $\tilde\mu=(\mu_2,\mu_3,\dots)$.  Then
\[
T_{w_\mu}y_{\mu'}\in \shif{\mu_1}{T_{w_{\tilde\mu}}y_{\tilde\mu'}}\hhh n.
\]
\end{propn}

To prove this, it will be helpful to generalise to \emph{skew partitions}.  For us, a skew partition is an ordered pair of partitions (written $\mu/\la$) where $\mu_i\gs\la_i$ for each $i$.  We define the Young diagram $[\mu/\la]$ to be the set difference $[\mu]\setminus[\la]$.  We identify two skew-partitions which have the same Young diagram, so that for example $(3,2^2)/(3,1^1)=(2^3\!,1)/(2,1^3)$.  If $\mu/\la$ is a skew-partition, then we define its conjugate $(\mu/\la)'$ to be $\mu'/\la'$.  We say that a skew partition $\tilde\mu$ is obtained from a skew partition $\mu$ \emph{by translation} if for some $i,j$ we have $[\tilde\mu]=\lset{(r+i,c+j)}{(r,c)\in[\mu]}$.  For most purposes $\tilde\mu$ and $\mu$ can be treated as equal in this case.

If $\mu$ is a skew partition, then the tableaux $\fkt^\mu$, $\fkt_\mu$, and hence the permutation $w_\mu$, can be defined in the obvious way, and so can $y_{\mu'}$.  Note that if $\tilde\mu$ is a skew-partition obtained from $\mu$ by translation, then $w_{\tilde\mu}=w_\mu$ and $y_{\tilde\mu'}=y_{\mu'}$.

Now Proposition \ref{tytop} follows by applying the next lemma $\mu_1$ times.

\begin{lemma}\label{tytopskew}
Suppose $\mu$ is a skew partition of $n>0$, and let $\nu$ be the skew partition obtained by removing the first node from the first non-empty row of $[\mu]$. Then $T_{w_\mu}y_{\mu'}\in \shif1{T_{w_\nu}y_{\nu'}}\hhh n$.
\end{lemma}

\begin{pf}
Suppose the node removed from $[\mu]$ to obtain $[\nu]$ is $(r,c)$.  Then $\fkt^\mu\shif1{w_\nu}$ is the standard tableau in which the number $1$ appears in position $(r,c)$, while the numbers $2,\dots,n$ appear in order down successive columns.  So, if the entry at the bottom of column $c{-}1$ is $m$, then $w_\mu=\shif1{w_\nu}w$, where $w$ is the cycle $(m\ m{-}1\ \dots\ 1)$.  Moreover, we have $l(w_\mu)=l(w_\nu)+l(w)$, so $T_{w_\mu}=\shif1{T_{w_\nu}}T_w=\shif1{T_{w_\nu}}T_1T_2\dots T_{m-1}$.

Now $y_{\mu'}$ is the commuting product of factors $y_k$, where $y_k=\sum_{s\in W(k)}(-q)^{-l(s)}T_s$ (with $W(k)$ being the group of permutations of the entries in column $k$ of $\fkt_\mu$).  $y_{\nu'}$ is the product of factors $\bar y_k$ defined similarly.  For $k>c$ we have $y_k=\shif1{\bar y_k}$, and $y_k$ commutes with $T_w$.  For $k=c$, we can factorise $y_k$ as $\shif1{\bar y_k}y$ using right coset representatives, and $\shif1{\bar y_k}$ again commutes with $T_w$.  For $k<c$ we have $y_k=\bar y_k$ and $T_wy_k=\shif1{y_k}T_w$.  And so
\[
T_{w_\mu}y_{\mu'}=\shif1{T_{w_\nu}}T_w\prod_ky_k=\shif1{T_{w_\nu}}\prod_k\shif1{\bar y_k}T_wy\in\shif1{T_{w_\nu}y_{\nu'}}\hhh n.\qhere{\sum_{k}}
\]
\end{pf}

\begin{eg}
Let $\mu=(3^2\!,2)/(1)$, so that $\nu=(3^2\!,2)/(2)$.  Then $w_\nu=(1\ 5\ 2)(4\ 6)$, and
\[
\fkt^\mu \shif1{w_\nu}=\gyoung(:;1;2,;3;4;5,;6;7)\,(2\ 6\ 3)(5\ 7)=\gyoung(:;1;6,;2;4;7,;3;5),\qquad w_\mu=\shif1{w_\nu}(1\ 3\ 2).
\]
We have $y_{\mu'}=y_1y_2y_3$ and $y_{\nu'}=\bar y_1\bar y_2\bar y_3$, where
\begin{alignat*}2
\bar y_1&=y_1=1-q^{-1}T_1,&\qquad\qquad&\\
\bar y_2&=1-q^{-1}T_3,\qquad &y_2&=\bar y_2\left(1-q^{-1}T_4+q^{-2}T_4T_3\right),\\
\bar y_3&=1-q^{-1}T_5,\qquad &y_3&=1-q^{-1}T_6=\shif1{\bar y_3}.
\end{alignat*}
So
\[
T_{w_\mu}y_{\mu'}=\shif1{T_{w_\nu}y_{\nu'}}T_1T_2\left(1-q^{-1}T_4+q^{-2}T_4T_3\right).
\]
\end{eg}

\begin{pf}[Proof of Proposition \ref{topmain}]
$S^\mu$ is generated by $x_\mu T_{w_\mu}y_{\mu'}$, so we need to show that
\[
\left(\sum_{i=1}^rc_ix_\mu\phi_{A_i}\right)T_{w_\mu}y_{\mu'}=0.
\]
Using Propositions \ref{rowtophom} and \ref{tytop}, this equals
\[
\left(\sum_{i=1}^rhc_i\shif{\mu_1}{x_{\tilde\mu}\phi_{\tilde A_i}}\right)\shif{\mu_1}{T_{w_{\tilde\mu}}y_{\tilde\mu'}}k,
\]
where $h,k\in\hhh n$.  Note that $h$ is independent of $i$: because $A_1,\dots,A_r$ all have the same first row, $\tilde A_1,\dots,\tilde A_r$ all have the same type.  So we have
\[
h\sh{\mu_1}\left(x_{\tilde\mu}\sum_{i=1}^rc_i\phi_{\tilde A_i}T_{w_{\tilde\mu}}y_{\tilde\mu'}\right)k.
\]
The term in the middle is $(x_{\tilde\mu}T_{w_{\tilde\mu}}y_{\tilde\mu'})\sum_{i=1}^rc_i\phi_{\tilde A_i}$, which is zero by assumption.
\end{pf}

\subs{Adding a row at the bottom}

We now consider a counterpart to Proposition \ref{topmain} in which we remove a row from the bottom of a tableau.

\begin{propn}\label{botmain}
Suppose $\mu\ptn n$ and $\la\compn n$, and let $l$ be maximal such that $\mu_l>0$.  Suppose $A_1,\dots,A_r$ are row-standard $\mu$-tableaux of type $\la$ which are identical in row $l$.  Write $\tilde\mu=(\mu_1,\dots,\mu_{l-1},0,0,\dots)$, and let $\tilde A_1,\dots,\tilde A_r$ be the $\tilde\mu$-tableaux obtained by deleting the $l$th row from each of $A_1,\dots,A_r$.  If $c_1,\dots,c_r\in\bbf$ are such that $\sum_{i=1}^rc_i\hat\phi_{\tilde A_i}=0$, then $\sum_{i=1}^rc_i\hat\phi_{A_i}=0$.
\end{propn}

This is proved in a very similar way to Proposition \ref{topmain}; the proof is in fact slightly simpler, because we do not need the functions $\sh s$, and we can more easily avoid skew partitions.  So we just give the two main propositions which are analogous to Propositions \ref{rowtophom} and \ref{tytop}, and leave the reader to fill in the remaining details.

\begin{propn}\label{rowbothom}
Suppose $A$ is a row-standard $\mu$-tableau of type $\la$, and let $l$ be maximal such that $\mu_l>0$.  Write $\tilde\mu=(\mu_1,\dots,\mu_{l-1},0,0,\dots)$, and let $\tilde A$ be the $\tilde\mu$-tableau obtained by deleting the last row of $A$.  Let $\tilde\la$ be the type of $\tilde A$.  Then we have
\[
x_\mu\phi_A=hx_{\tilde\mu}\phi_{\tilde A}
\]
where $h\in\hhh n$ depends only on $\mu,\la,\tilde\la$.
\end{propn}

\begin{propn}\label{tybot}
Suppose $\mu$ is a partition of $n$, and write $\tilde\mu=(\mu_1,\dots,\mu_l,0,0,\dots)$, where $l$ is maximal such that $\mu_l>0$.  Then
\[
T_{w_\mu}y_{\mu'}\in T_{w_{\tilde\mu}}y_{\tilde\mu'}\hhh n.
\]
\end{propn}

\section{Expressing a homomorphism as a linear combination of semistandard homomorphisms}\label{s:algo}

\subs{Two-part partitions}

We now explain how to use Theorem \ref{main} to write a row-standard homomorphism $\hat\phi_A$ as a linear combination of semistandard homomorphisms.  We begin by considering two-part partitions.

Suppose $\mu=(m,n-m)$ is a two-part partition, and $A$ is a row-standard $\mu$-tableau of arbitrary type.  If $A$ is not semistandard, then there is some $i\ls n-m$ such that the $(1,i)$-entry of $A$ is at least as large as the $(2,i)$-entry.  Suppose the $(1,i)$-entry equals $j$.  Now define the following multisets:
\begin{align*}
R&=\lset{k\in A^1}{k<j};\\
S&=\lset{k\in A^1}{k\gs j}\sqcup\lset{k\in A^2}{k\ls j};\\
T&=\lset{k\in A^2}{k>j}.
\end{align*}
Then $|R|<i$ and $|T|\ls n-m-i$, so $R,S,T$ satisfy the conditions of Theorem \ref{main}.  Moreover, the tableau $A$ arises as $A[U_0,V_0]$, where
\[
U_0=\lset{k\in A^1}{k\gs j},\qquad V_0=\lset{k\in A^2}{k\ls j};
\]
note in particular that $U_0$ consists of the $|U_0|$ largest elements of $S$.  Note also that all the elements of $R$ are strictly less than all the elements of $U_0$, and the same is true for $V_0$ and $T$, so the coefficient of $\hat\phi_A$ in Theorem \ref{main} is $1$.

So Theorem \ref{main} allows us to express $\hat\phi_A$ as a linear combination of homomorphisms $\hat\phi_B$, for row-standard tableaux $B$ satisfying
\[
\sum_{k\in B^1}k<\sum_{k\in A^1}k.
\]
We can repeat the process for each $B$ which is not semistandard, and eventually express $\hat\phi_A$ as a linear combination of semistandard homomorphisms. The above inequality guarantees that this process will terminate.

\begin{eg}
For this example, we abuse notation and write a homomorphism $\hat\phi_B$ just as $B$.  Suppose
\Yvcentermath1
\[
A=\young(12234,1333).
\]
Applying Theorem \ref{main} with
\[
R=\emptyset,\qquad S=\{1,1,2,2,3,4\},\qquad T=\{3,3,3\},
\]
we get
\[
\hat\phi_A=-[4]\young(11224,3333)-\young(11234,2333)-q^3\young(11223,3334).
\]
Applying Theorem \ref{main} again, we have
\[
\young(11234,2333)=-[2]\young(11223,3334)-[2]\young(11224,3333)-\young(11233,2334)
\]
so that
\[
\hat\phi_A=-(q^2+q^3)\young(11224,3333)+(1+q-q^3)\young(11223,3334)+\young(11233,2334).
\]
\end{eg}

\subs{Arbitrary partitions}

Now we consider generalising this to arbitrary partitions, using the results of Section \ref{s:rowrem}.  By applying Propositions \ref{topmain} and \ref{botmain} repeatedly, we obtain the following theorem.

\begin{thm}\label{mainrowrem}
Suppose $\mu\ptn n$ and $\la\compn n$, and $1\ls l<m$.  Let $\tilde\mu$ be the partition $(\mu_l,\mu_{l+1},\dots,\mu_m,0,0,\dots)$.  Suppose $A_1,\dots,A_r$ are row-standard $\mu$-tableaux of type $\la$ which are identical outside rows $l,l{+}1,\dots,m$, and let $\tilde A_1,\dots,\tilde A_r$ be the $\mu$-tableaux obtained by deleting rows $1,\dots,l{-}1$ and $m{+}1,m{+}2,\dots$ of $A_1,\dots,A_r$.  If $c_1,\dots,c_r\in\bbf$ are such that $\sum_{i=1}^rc_i\hat\phi_{\tilde A_i}=0$, then $\sum_{i=1}^rc_i\hat\phi_{A_i}=0$.
\end{thm}

Now suppose $A$ is a row-standard $\mu$-tableau.  If $A$ is not semistandard, then there must be some $l$ such that the tableau $\tilde A$ consisting of rows $l,l{+}1$ of $A$ is not semistandard.  Using the algorithm described above, we can express $\hat\phi_{\tilde A}$ as a linear combination of homomorphisms of the form $\hat\phi_{\tilde B}$ where each $\tilde B$ satisfies $\sum_{k\in \tilde B^1}k<\sum_{k\in\tilde A^1}k$.  Applying Theorem \ref{mainrowrem}, we can express $\hat\phi_A$ as a linear combination of homomorphisms $\hat\phi_B$ where each $B$ satisfies
\[
\sum_{i\gs1}\sum_{k\in B^i}ik>\sum_{i\gs1}\sum_{k\in A^i}ik.
\]
So by induction we can express $\phi_A$ as a linear combination of semistandard homomorphisms.

\subs{Computing the space of homomorphisms between two Specht modules}

We now explain briefly how to compute the space of $\hhh n$-homomorphisms between two Specht modules when $q\neq-1$, which is the main motivation for proving the results in this paper.  A more detailed discussion of this topic can be found in the paper \cite{L3} by Lyle.

Suppose $\mu,\la\ptn n$. Theorem \ref{semi}(2) says that the set
\[
\lset{\hat\phi_A}{A\text{ a semistandard $\mu$-tableau of type $\la$}}
\]
is a basis for $\hom_{\hhh n}(S^\mu,M^\la)$ as long as $q\neq-1$.  Since $S^\la\subseteq M^\la$, an $\hhh n$-homomorphism from $S^\mu$ to $S^\la$ is simply a linear combination of semistandard homomorphisms whose image lies inside the Specht module $S^\la$.  Fortunately, there is a straightforward way to check this condition: $S^\la$ may be characterised \cite[Theorem~7.5]{dj} as the intersection of the kernels of certain tableau homomorphisms $\psi_{i,j}$ from $M^\la$ to other permutation modules $M^\nu$.  Given a tableau homomorphism $\phi_A:M^\mu\to M^\la$, it is known how to compose $\phi_A$ with $\psi_{i,j}$; this result is \cite[Proposition~2.14]{sl}, generalising the result \cite[Lemma 5]{fm} proved for the symmetric group.  However, this composition is given as a linear combination of tableau homomorphisms which are not always semistandard (even if $A$ is semistandard).  Since the set $\lset{\hat\phi_A}{A\in\caltr\mu\la}$ is linearly dependent in general, this can make it difficult to tell whether a given linear combination $\sum_Ac_A\hat\phi_A\psi_{i,j}$ is non-zero.  But using the results in this paper, one can quickly express this homomorphism as a linear combination of semistandard homomorphisms, and hence discover whether it is non-zero.

So one has a reasonably fast algorithm for computing the space of $\hhh n$-homomorphisms between two Specht modules when $q\neq-1$.  The author has implemented this in a small package for the GAP4 system \cite{gap4}.  This package has already yielded new results on homomorphisms between Specht modules \cite{allreds,sl2d,df}.

\end{document}